\documentclass{amsart}
\usepackage{amssymb}

\usepackage{graphicx}
\usepackage{amscd}

\newtheorem{theorem}{Theorem}
\theoremstyle{plain}

\newtheorem{corollary}{Corollary}

\newtheorem{lemma}{Lemma}

\newtheorem{proposition}{Proposition}
\newtheorem{remark}{Remark}

\numberwithin{equation}{section}

\begin{document}
\title[mapping class group]{The mapping class group of a disk with infinitely many holes }
\author{Paul Fabel}
\address{Drawer MA\\
Department of Mathematics and Statistics\\
Mississippi State, Ms. 39762}
\email{fabel@ra.msstate.edu}
\urladdr{http://www2.msstate.edu/\symbol{126}fabel/}
\date{March 4, 2003}

\begin{abstract}
A left orderable completely metrizable topological group is exhibited
containing Artin's braid group on infinitely many strands. The group is the
mapping class group (rel boundary) of the closed unit disk with a sequence
of interior punctures converging to the boundary. This resolves an issue
suggested by work of Dehornoy.
\end{abstract}

\maketitle

\section{Introduction}

Artin's braid group $B_{n}$ is the fundamental group of the configuration
space consisting of planar sets with precisely $n$ elements.

Artin's braid group on infinitely many strands $B_{\infty }$ is the direct
limit of $B_{n}$ under monomorphisms $k_{n}:B_{n}\hookrightarrow B_{n+1}$
attaching a `trivial strand'. Each element of $B_{\infty }$ can be seen as a
countable collection of disjoint arcs in $[0,1]\times R^{2}$ connecting $%
\{0\}\times \{(1,0),(2,0),...\}$ to $\{1\}\times \{(1,0),(2,0),..\}$ such
that all but finitely many of the arcs are line segments.

The group $B_{\infty }$ admits two natural incomplete topologies yielding
distinct topological groups each of which is homeomorphic to the rational
numbers.

In \cite{fab} a completely metrizable group $\overline{B_{\infty }}$ is
constructed completing $B_{\infty }$ with the following topology: A sequence
of braids $b_{n}\in B_{\infty }$ converges to $id\in B_{\infty }$ iff for
each $M$ there exists $N$ such that if $n\geq N$ then the first $M$ strands
of $b_{n}$ determine the trivial braid on $M$ strands. A typical element of
the completion is a `wild braid' and is not determined by an ambient isotopy
of the plane. The pure subgroup is precisely the inverse limit of Artin's
pure braid groups.

However there is another natural topology on $B_{\infty }$ determined by a
left ordering explored extensively in \cite{deh} \cite{funk}, in which $%
b_{1},b_{2},...$ converges to $id$ iff for each $M$ there exists $N$ such
that for $n\geq N$ the first $M$ strands of $b_{n}$ determine the trivial
braid on $M$ strands \textbf{and} are unlinked from the remaining strands of 
$b_{n}.$

In reference to $B_{\infty }$ with this order topology, Dehornoy \cite{deh}
has commented ``... a group structure on the completion,...would be more
satisfying than the monoid structure of $EB_{\infty }.$''

To this end, we exhibit a left ordered completely metrizable mapping class
group $M(H_{\infty })$ in which $B_{\infty }$ with the order topology from 
\cite{deh} is a densely embedded subgroup. Choosing a sequence $z_{n}\in
int(D^{2})$ converging to a boundary point of a closed disk $D^{2}\subset
R^{2},$ $H_{\infty }$ is the group of homeomorphisms of $D^{2}$ which fix $%
\partial D^{2}$ pointwise and which leave the set $\{z_{1},z_{2},...\}$
invariant.

In contrast to $\overline{B}_{\infty }$ each element of $M(H_{\infty })$ is
a `tame braid' determined by an isotopy of the disk. The natural map $%
j:M(H_{\infty })\hookrightarrow \overline{B_{\infty }}$ is a continuous
monomorphism but not an embedding. Both $M(H_{\infty })$ and $\overline{B}%
_{\infty }$ are homeomorphic to the irrational numbers.

Complications arise from the fact that $H_{\infty }$ is not locally
connected. In particular the path component of the identity $H_{\infty }^{0}$
is not an open subgroup of $H_{\infty }$ and some care is required to verify
(Lemma \ref{close}) that $H_{\infty }^{0}$ is in fact closed in $H_{\infty
}. $

A general procedure is developed (Corollary \ref{cor2}) for generating a
metrically complete group given a normal subgroup $H$ of a subgroup $G$ of
the autohomeomorphism group of a metric compactum.

This determines the following complete metric $D$ on $M(H_{\infty }):$

Define $D(gH_{\infty }^{0},fH_{\infty }^{0})=D_{H}(gH_{\infty
}^{0},fH_{\infty }^{0})+D_{H}(g^{-1}H_{\infty }^{0},f^{-1}H_{\infty }^{0})$
where $D_{H}$ is the (incomplete) Hausdorf metric determined by the condition

$D_{H}(pH_{\infty }^{0},qH_{\infty }^{0})<\varepsilon $ iff there exists $p^{%
\symbol{94}}\in pH_{\infty }^{0}$ and $q^{\symbol{94}}\in qH_{\infty }^{0}$
such that $\left| p^{\symbol{94}}(z)-q^{\symbol{94}}(z)\right| <\varepsilon
\forall z\in D^{2}.$

Dehornoy's left ordering on $B_{n}$ induces a left ordering on $M(H_{\infty
}).$ It is verified in Theorem \ref{main} that the order topology is
compatible with the quotient topology of $M(H_{\infty }).$

\section{Basic definitions and notation}

Throughout this paper \bigskip all function spaces will have the compact
open topology and all identification spaces will have the quotient topology. 
$C(A,B)$ denotes the maps from $A$ to $B,$ and $H(A,A)$ denotes the
homeomorphisms from $A$ onto $A.$

If $G$ is a topological group then $G^{0}\subset G$ denotes the path
component of the identity in $G$, and $MG$ denotes, mapping class group of $%
G,$ the quotient space $G/G^{0}.$

Let $D^{2}\subset R^{2}$ denote the closed disk of radius $1$ centered at $%
(1,0).$

Let $Z_{n}=\{(1,0),(\frac{1}{2},0),..(\frac{1}{n},0)\}\subset int(D^{2}).$

Let $Z=\cup _{n=1}^{\infty }Z_{n}$

Let $\alpha _{n}=\{(x,y)\in D^{2}|x=\frac{\frac{1}{n+1}+\frac{1}{n}}{2}\}.$
Thus $\alpha _{n}$ is a vertical line segment halfway between $z_{n+1}$ and $%
z_{n}.$

Let $D_{n}=\{(x,y)\in D^{2}|\frac{\frac{1}{n+1}+\frac{1}{n}}{2}\leq x\}.$

Given $n>k$ let $D_{n,k}=\overline{D_{n}\backslash D_{k}}.$

Let $H_{\infty }=\{h\in H(D^{2},D^{2})|h(Z)=Z$ and $h_{\partial
D^{2}}=id_{\partial D^{2}}\}.$

Let $H_{n}=\{h\in H(D^{2},D^{2})|h(Z_{n})=Z_{n}$ and $h_{\partial
D^{2}}=id_{\partial D^{2}}\}.$

Let $H_{n}^{\symbol{94}}=\{h\in H_{n}|h_{D^{2}\backslash
D_{n}}=id_{D^{2}\backslash D_{n}}\}$

\section{Artin's braid groups}

A braid $b_{n}\in B_{n}$ can be seen as the restriction to $n$ vertical line
segments of a level preserving homeomorphism of the $3$-cell $D^{2}\times
\lbrack 0,1].$ Formally we select from several equivalent definitions (see 
\cite{r1} for a quick tour) and define $B_{n},$ Artin's braid group on $n$
strands to be the mapping class group $M(H_{n}).$ The following propositions
are well known. See \cite{Birm} for more details.

\begin{proposition}
\label{prop1}The inclusion map $i_{n}:H_{n}^{\symbol{94}}\rightarrow H_{n}$
induces an isomorphism $i_{n}^{\ast }:M(H_{n}^{\symbol{94}})\rightarrow
M(H_{n}).$
\end{proposition}

\begin{proposition}
\label{prop2}The inclusion map $j_{n}:H_{n}^{\symbol{94}}\rightarrow
H_{n+1}^{\symbol{94}}$ induces a monomorphism $j_{n}^{\ast }:M(H_{n}^{%
\symbol{94}})\hookrightarrow M(H_{n+1}^{\symbol{94}})$
\end{proposition}

Now define Artin's braid group on infinitely many strands $B_{\infty }$ as
the direct limit of $B_{n}$ taken over the monomorphisms $%
k_{n}=i_{n+1}^{\ast }j_{n}^{\ast }i_{n}^{\ast -1}:B_{n}\hookrightarrow
B_{n+1}$

\begin{lemma}
The monomorphisms $i_{n+1}^{\ast }j_{n}^{\ast }i_{n}^{\ast
-1}:B_{n}\rightarrow B_{n+1}$ induce a monomorphism $i_{\infty }:B_{\infty
}\hookrightarrow M(H_{\infty })$.
\end{lemma}

\begin{proof}
By proposition \ref{prop2} each element of $\lim_{\rightarrow }M(H_{n})$ can
be represented in the form $gH_{n}^{0},gH_{n+1}^{0},...$ with $g\in H_{n}^{%
\symbol{94}}$ and $gH_{n}^{0}\notin im(k_{n-1}).$

Thus $i_{\infty }(gH_{n}^{0},gH_{n+1}^{0},)=gH_{\infty }^{0}.$ Suppose $%
i_{\infty }(gH_{n}^{0},gH_{n+1}^{0},)=$ $gH_{\infty }^{0}=H_{\infty }^{0}.$%
Then $g\in H_{\infty }^{0}\subset \cap _{m=1}^{\infty }H_{m}^{0}$.

Thus $n=1$ and $(gH_{n}^{0},gH_{n+1}^{0},...)=(H_{1}^{0},H_{2}^{0},....).$
Hence $i_{\infty }$ is one to one.
\end{proof}

\section{Building isotopies in $H_{\infty }$}

This section produces Corollaries \ref{nicecor} and \ref{close}.

\begin{remark}[Alexander trick]
Any two homeomorphisms of $D^{2}$ which agree on $\partial D^{2}$ and a
point $p\in int(D^{2})$ are isotopic rel $\partial D^{2}\cup p.$
\end{remark}

\begin{remark}
\label{ext}Suppose $\alpha \subset D^{2}\backslash Z_{n}$ is a closed arc $.$
Suppose $g_{t}:\alpha \hookrightarrow D^{2}\backslash Z_{n}$ is an isotopy
such that $g_{0}=id_{\alpha },\partial (im(g_{t}))\subset \partial D^{2}$
and $int(im(g_{t}))\subset int(D^{2}\backslash Z_{n}).$ Then there exists an
isotopy $G_{t}=D^{2}\rightarrow D^{2}$ such that $G_{0}=id_{D^{2}},(G_{t})_{%
\alpha }=g_{t}$ and $G_{t}\in H_{n}^{0}.$
\end{remark}

\begin{lemma}
\label{kill}Suppose $g\in H_{\infty },$ $g_{D_{i}}=id_{D_{i}}$ and $\forall
n $ $id_{\alpha _{i+1}}$ is isotopic ( rel $\partial \alpha _{i+1}$) in $%
D^{2}\backslash Z_{n}$ to $g_{\alpha _{i+1}}.$ Then there exists $f\in
gH_{\infty }^{0}$ such that $f_{D_{i+1}}=id_{D_{i+1}}.$
\end{lemma}

\begin{proof}
Choose $n>i+1$ such that $g(\alpha _{i+1})\cap \alpha _{n}=\emptyset .$ Let $%
\beta _{t}$ be an isotopy rel $\partial \alpha _{i+1}$ in $D^{2}\backslash
Z_{n}$ such that $\beta _{0}=id_{\alpha _{i+1}}$ and $\beta _{1}=g_{\alpha
_{i+1}}.$ By Remark \ref{ext} let $h_{t}$ be an isotopy in $H_{n}^{0}$
extending $\beta _{t}$ such that $h=id_{D^{2}}$. Let $\phi :D^{2}\rightarrow
D_{n,i}$ be a homeomorphism fixing pointwise $\alpha _{i+1}\cup
(Z_{n}\backslash Z_{i}).$ Let $g_{t}^{\ast }=id_{(D^{2}\backslash
D_{n})}\cup \phi h_{t}\phi ^{-1}\cup id_{(D_{k})}$. Note $g_{t}\in H_{\infty
}^{0}$ and $g_{t}$ fixes $D_{i}$ pointwise. Now consider $%
p_{t}=g(g_{t}^{\ast })^{-1}.$ Note $p_{0}=g$ and $p_{1}$ fixes $\alpha
_{i+1}\cup \alpha _{i}$ pointwise. Hence by the Alexander trick we may
isotop $p_{1}$ within $H_{\infty }^{0}$ to a map $f$ such that $%
f_{D_{i+1}}=id_{D_{i+1}}.$
\end{proof}

\begin{corollary}
\label{nicecor} Suppose $g\in H_{\infty }$ and for each $i\leq k$ and for
each $n$ $\alpha _{i}\in (gH_{n}^{0})(\alpha _{i}).$ Then there exists $g^{%
\symbol{94}}\in gH_{\infty }^{0}$ such that $g_{D_{k}}^{\symbol{94}%
}=id_{D_{k}}.$
\end{corollary}

\begin{proof}
Successively apply Lemma \ref{kill} $k$ times to construct a map $\psi
:[0,k]\rightarrow gH_{\infty }^{0}$ satisfying $\psi (t)_{D_{i}}=id_{D_{i}}$
whenever $t\leq i$ and $i\in \{1,2,...k\}.$ Let $g^{\symbol{94}}=\psi (k).$
\end{proof}

\begin{lemma}
\label{exiso}Suppose $h\in \cap H_{n}^{0}$ and $h_{D_{k}}=id_{D_{k}}.$ Then
there exists an isotopy $h_{t}$ in $\cap H_{n}^{0}$ such that $h_{0}=h,$ $%
h_{t(D_{k})}=id_{(D_{k})}$ and $h_{1(D_{k+1})}=id_{(D_{k+1})}.$
\end{lemma}

\begin{proof}
Choose $n>k+1$ such that $\alpha _{n}\cap h(\alpha _{k+1})=\emptyset .$ Let $%
g_{t}$ be a isotopy in $H_{n}^{0}$ such that $id=g_{0}$ and $g_{1}=h$. Let $%
\phi :D^{2}\rightarrow \overline{D_{n}\backslash D_{k}}$ be a homeomorphism
such that $\phi _{\alpha _{k+1}\cup (Z_{n}\backslash Z_{k})}=id_{\alpha
_{k+1}\cup (Z_{n}\backslash Z_{k})}.$ Let $h_{t}^{\ast
}=id_{(D^{2}\backslash D_{n})}\cup \phi g_{t}\phi ^{-1}\cup id_{D_{k}}.$ Let 
$p_{t}=h(h_{t}^{\ast })^{-1}.$ Note $p_{1}$ fixes $\alpha _{k}\cup \alpha
_{k+1}$ pointwise. Now isotop $p_{1(D_{k+1},k)}$ to $id_{D_{k+1},k}$ via the
Alexander trick while leaving $D_{k}\cup (D^{2}\backslash D_{k+1})$ alone.
\end{proof}

\begin{corollary}
\label{close}$H_{\infty }^{0}$ is closed in $H_{\infty }$ and $\cap
_{n=1}^{\infty }H_{n}^{0}=H_{\infty }^{0}.$
\end{corollary}

Suppose $h\in \cap H_{n}^{0}.$ By repeated application of Lemma \ref{exiso}
construct a map $\psi :[0,\infty )\rightarrow \cap H_{n}^{0}$ satisfying $%
\psi (0)=h$ and $\psi (t)_{D_{n}}=id_{D_{n}}$ for $t\leq n.$ Thus $\psi $
can be continuously extended to $[0,\infty ]$ with $\psi (\infty )=id.$
Hence $\cap H_{n}^{0}\subset H_{\infty }^{0}.$ Conversely $H_{\infty
}\subset H_{n}\forall n$ and hence $H_{\infty }^{0}\subset H_{n}^{0}\forall
n.$ Hence $H_{\infty }^{0}=\cap _{n=1}^{\infty }H_{n}^{0}.$ Moreover $%
H_{\infty }^{0}$ is closed in $H_{\infty }$ since its the intersection of
closed sets $H_{n}^{0}\subset H(D^{2},D^{2}).$

\section{Recognizing complete quotient groups.}

Expanding an elementary argument ( Proposition 1.3.10 \cite{vanmill}) that
the autohomeomorphism $H(X,X)$ of a compact metric $X$ is completely
metrizable, given a subgroup $G\subset H(X,X)$ with normal subgroup $%
H\vartriangleleft G,$ we construct a completely metrizable group $(\overline{%
G}/\overline{H})^{\symbol{94}}$ such that if $G$ and $H$ are closed in $%
H(X,X)$ then $(\overline{G}/\overline{H})^{\symbol{94}}\simeq G/H.$ The
proof of Theorem \ref{comgrp} makes use of the Baire category theorem, in a
complete metric space the nested intersection of countably many open dense
sets is dense, and in particular nonempty.

\begin{remark}
Suppose $(X,d)$ is a compact metric space. Fixing $\varepsilon >0$ $,$ a map 
$f\in C(X,X)$ is an $\mathbf{\varepsilon }$ \textbf{map} if $%
diam(f^{-1}(y))<\varepsilon $ $\forall y\in X.$ The $\varepsilon $ maps form
an open subspace of $C(X,X).$
\end{remark}

\begin{lemma}
\label{conv}Suppose $X$ is a compact metric space and $%
\{f,f_{1}f_{2},..,p_{1},p_{2},.....\}\subset C(X,X).$ Suppose $%
fp_{n}\rightarrow id$ and $p_{n}f_{n}\rightarrow id.$ Then $f_{n}\rightarrow
f.$
\end{lemma}

\begin{proof}
Given $g\in C(X,X)$ denote the \textbf{graph} of $g$ as $gr(g)=\{(x,y)\in
X\times X|y=g(x)\}.$ Let $f_{n}^{\ast }$ be a subsequence of $f_{n}$ such
that $gr(f_{n}^{\ast })$ converges (with respect to the Hausdorf metric over
compacta in $X\times X)$ to a compactum $\mathcal{X}\subset X\times X.$ Let $%
p_{n}^{\ast }$ denote the corresponding subsequence of $p_{n}.$ ( Thus if $%
f_{n}^{\ast }=f_{n_{m}}$ then $p_{n}^{\ast }=p_{n_{m}}$). Suppose $b=f(a).$
Let $a=p_{n}^{\ast }f_{n}^{\ast }(x_{n}).$ Then $x_{n}\rightarrow a.$
Moreover $f(a)=\lim f(p_{n}^{\ast }f_{n}^{\ast }(x_{n}))=\lim (fp_{n}^{\ast
})(f_{n}^{\ast }(x_{n}))=\lim f_{n}^{\ast }(x_{n}).$ This shows $%
gr(f)\subset \mathcal{X}.$ Suppose $(a,b)\in \mathcal{X}.$ Choose a
subsequence $f_{n}^{\symbol{94}}$ of $f_{n}^{\ast }$ such that $b=\lim
f_{n}^{\symbol{94}}(a).$ Let $p_{n}^{\symbol{94}}$ denote the corresponding
subsequence of $p_{n}.$ Then $b=\lim f_{n}^{\symbol{94}}(a)=\lim (fp_{n}^{%
\symbol{94}})(f_{n}^{\symbol{94}}(a))=\lim f((p_{n}^{\symbol{94}}f_{n}^{%
\symbol{94}}(a))=\lim f(a)=f(a).$ This shows every subsequential limit $%
\mathcal{X}$ of $gr(f_{n})$ satisfies $\mathcal{X}=gr(f).$ Hence $%
f_{n}\rightarrow f.$
\end{proof}

\begin{lemma}
\label{compat}Suppose $G$ is a topological group with compatible metric $d$.
Then the function $D:G\times G\rightarrow G$ defined via $%
D(f,g)=d(f,g)+d(f^{-1},g^{1})$ is a compatible metric on $G$.
\end{lemma}

\begin{proof}
Note $D(f,g)=0$ iff $d(f,g)+d(f^{-1},g^{-1})=0$ iff $f=g$ (since $d\geq 0$).
Symmetry is immediate. $%
D(f,g)+D(g,p)=d(f,g)+d(g,p)+d(f^{-1},g^{-1})+d(g^{-1},p^{-1})\geq
d(f,p)+d(f^{-1},p^{-1})=D(f,p).$ To prove $D$ determines the same topology
we will show that $(G,d)$ and $(G,d)$ have the same convergent sequences. If 
$D(g_{n},g)\rightarrow 0$ then immediately $d(g_{n},g)\rightarrow 0$ since $%
d\leq D.$ Suppose $d(g_{n},g)\rightarrow 0.$ Then $d(g_{n}^{-1},g^{-1})%
\rightarrow 0$ since inversion is continuous in $G.$ Thus $%
D(g,g_{n})\rightarrow 0.$
\end{proof}

\begin{lemma}
\label{lem1}Suppose $h:X\rightarrow X$ is a homeomorphism of the compact
metric space $(X,d).$ Then the map $\phi _{h}:C(X,X)\rightarrow C(X,X)$
defined via $\phi _{h}(f)=fh$ is an isometry with the uniform metric $D$ on $%
C(X,X).$
\end{lemma}

\begin{proof}
Suppose $f\in C(X,X).$ Then $\phi _{h}(fh^{-1})=f.$ Thus $\phi _{h}$ is
surjective. Suppose $f,g\in C(X,X)$ and $D(f,g)=d(f(x),g(x)).$ Let $%
y=h^{-1}(x).$ Then $D(fh,gh)\geq d(fh(y),d(gh(y)))=D(f,g).$

Conversely suppose $D(fh,gh)=d(fh(y),gh(y)).$ Then, letting $x=h(y),$ $%
D(f,g)\geq d(f(x),g(x))=D(fh,gh).$
\end{proof}

\begin{lemma}
\label{lem2}Suppose $H$ a group of isometries of the metric space $(Y,d).$
Then sets of the form $\overline{H(y)}$ determine a partition of $Y.$ The
function $\Pi :Y\rightarrow 2^{Y}$ defined via $\Pi (y)=\overline{H(y)}$ is
a quotient map with respect to the Hausdorf metric $D$. Moreover $D(%
\overline{H(y)},\overline{H(x)})<\varepsilon $ iff there exists $x^{\symbol{%
94}}\in \overline{H(x)},$ and $y^{\symbol{94}}\in \overline{H(y)}$ such that 
$d(x^{\symbol{94}},y^{\symbol{94}})<\varepsilon .$
\end{lemma}

\begin{proof}
Suppose $x^{\symbol{94}}\in \overline{H(x)}$ $y^{\symbol{94}}\in \overline{%
H(y)}$ and $d(x^{\symbol{94}},y^{\symbol{94}})<\varepsilon .$ Choose $%
\varepsilon ^{\symbol{94}}$ such that $d(x^{\symbol{94}},y^{\symbol{94}%
})<\varepsilon ^{\symbol{94}}<\varepsilon .$ Suppose $z\in \overline{H(x)}.$
Suppose $\delta >0.$ Choose $\{g,h,p\}\subset H$ such that $%
d(g(x),h(y))<\varepsilon ^{\symbol{94}}$ and $d(z,p(x))<\delta $. Then $%
d(z,pg^{-1}h(y))\leq d(z,p(x))+d(p(x),pg^{-1}h(y))<\delta
+d(x,g^{-1}h(y))=\delta +d(g(x),h(y))<\delta +\varepsilon ^{\symbol{94}}.$
Since $\delta $ was arbitrary this shows each point of $\overline{H(x)}$ is
no further than $\varepsilon ^{\symbol{94}}$ from some point of $\overline{%
H(y)}.$ By a symmetric argument we conclude $D(\overline{H(y)},\overline{H(x)%
})\leq \varepsilon ^{\symbol{94}}<\varepsilon .$ Thus $D(\overline{H(y)},%
\overline{H(x)})<\varepsilon $ iff there exists $x^{\symbol{94}}\in 
\overline{H(x)},$ and $y^{\symbol{94}}\in \overline{H(y)}$ such that $d(x^{%
\symbol{94}},y^{\symbol{94}})<\varepsilon .$ In particular if $z\in 
\overline{H(y)}\cap \overline{H(x)}$ then $D(\overline{H(y)},\overline{H(x)}%
)=0$ and hence $\overline{H(x)}=\overline{H(y)}.$ Thus sets of the form $%
\overline{H(x)}$ determine a partition of $X.$

Finally, to prove $\Pi $ is a quotient map first observe that $\Pi $ is
continuous since $D(\Pi (x),\Pi (y))\leq d(x,y).$ Next suppose $\Pi ^{-1}(A)$
is closed.

Suppose $D(\overline{H(x)},\overline{H(x_{n})})\rightarrow 0$ with $%
\overline{H(x_{n})}\in A.$ Choose $y_{n}\in \overline{H(x_{n})}$ with $%
d(x,y_{n})\rightarrow 0.$ Thus $x\in \Pi ^{-1}(A).$ Hence $\Pi (x)=\overline{%
H(x)}\in A.$ Thus $A$ is closed.
\end{proof}

Combining the previous two Lemmas shows that any group of homeomorphisms $H$
acting on a compact metric space $X$ generates a reasonable quotient of $%
C(X,X)$ denoted $C_{H}(X,X)$ as follows.

\begin{corollary}
\label{cor1}Suppose $H$ is a group of homeomorphisms of the compact metric
space $(X,d)$. Then, for each $f\in C(X,X),$the sets of the form $\overline{%
f(H)}$ determine a partition of $C(X,X)$. If $C_{H}(X,X)$ denotes the
associated quotient space then the Hausdorf metric is compatible with the
quotient topology.
\end{corollary}

\begin{proof}
By Lemma \ref{lem1} $H$ acts isometrically on $C(X,X).$ Now apply Lemma \ref
{lem2}.
\end{proof}

Throughout the remainder of this section $G$ is a group of homeomorphisms of
the compact metric space $X,$ $H$ is a normal subgroup of $G$ and $\overline{%
G}$ denotes the closure of $G$ in $C(X,X).$ $\overline{G}/\overline{H}$
denotes the subspace of $C_{H}(X,X)$ consisting of elements of the form $%
\overline{fH}$ with $f\in \overline{G}.$

\begin{lemma}
\label{th1} Define $\ast :\overline{G}/\overline{H}\times \overline{G}/%
\overline{H}\rightarrow \overline{G}/\overline{H}$ via $\ast ((\overline{fH}%
),(\overline{pH}))=\overline{fpH}.$ Then $\ast $ is well defined,
continuous, and determines that $\overline{G}/\overline{H}$ is a \textbf{%
monoid. }
\end{lemma}

\begin{proof}
To prove $\ast $ is well defined suppose $f^{\symbol{94}}\in (\overline{fH})$
and $p^{\symbol{94}}\in (\overline{pH}).$ It suffices to show $f^{\symbol{94}%
}p^{\symbol{94}}\in \overline{fpH}.$ Choose $\{f_{n},p_{n}\}\in G$ with $%
f_{n}\rightarrow f$ and $p_{n}\rightarrow p.$ Thus $\overline{f_{n}H}%
\rightarrow \overline{fH}$ and $\overline{p_{n}H}\rightarrow \overline{pH}.$
Since $f_{n}H$ is dense in $\overline{f_{n}H}$ we may choose $f_{n}^{\symbol{%
94}}\in f_{n}H$ and $p_{n}^{\symbol{94}}\in p_{n}H$ with $f_{n}^{\symbol{94}%
}\rightarrow f^{\symbol{94}}$ and $p_{n}^{\symbol{94}}\rightarrow p^{\symbol{%
94}}.$ Let $f_{n}^{\symbol{94}}=f_{n}h_{n}$ and $p_{n}^{\symbol{94}%
}=p_{n}h^{n}$ with $\{h_{n},h^{n}\}\in H.$ Then,since $H$ is normal $f_{n}^{%
\symbol{94}}p_{n}^{\symbol{94}%
}=(f_{n}h_{n})p_{n}h^{n}=f_{n}p_{n}(p_{n}^{-1}h_{n}p_{n})h^{n}\in
f_{n}p_{n}H.$ Note $\overline{(f_{n}p_{n})H}\rightarrow \overline{fpH}$
since $f_{n}p_{n}\rightarrow fp.$ Thus $f^{\symbol{94}}p^{\symbol{94}}=\lim
f_{n}^{\symbol{94}}p_{n}^{\symbol{94}}\in $ $\overline{fpH}.$ Hence $\ast $
is well defined. To prove $\ast $ is continuous suppose $\overline{f_{n}H}%
\rightarrow \overline{fH}$ and $\overline{p_{n}H}\rightarrow \overline{pH}.$
Choose $f_{n}^{\symbol{94}}\in \overline{f_{n}H}$ and $p_{n}^{\symbol{94}%
}\in \overline{p_{n}H}$ such that $f_{n}^{\symbol{94}}\rightarrow f$ and $%
p_{n}^{\symbol{94}}\rightarrow p.$ Thus $\overline{f_{n}^{\symbol{94}}p_{n}^{%
\symbol{94}}H}\rightarrow \overline{fpH}.$ Note $\overline{H}\ast \overline{%
fH}=\overline{H}\overline{fH}=\overline{fH}\ast \overline{H}.$
\end{proof}

\begin{corollary}
\label{cor2}Let $(\overline{G}/\overline{H})^{\symbol{94}}$denote the
subspace of $\overline{G}/\overline{H}$ consisting of invertible elements. (
i.e. $\overline{fH}\in (\overline{G}/\overline{H})^{\symbol{94}}$ iff there
exists $\overline{pH}\in \overline{G}/\overline{H}$ such that $\overline{pfH}%
=\overline{fpH}=\overline{H}$) . Then $(\overline{G}/\overline{H})^{\symbol{%
94}}$ is a completely metrizable topological group.\ 
\end{corollary}

\begin{proof}
$(\overline{G}/\overline{H})^{\symbol{94}}$ is a group since $\overline{G}/%
\overline{H}$ is a monoid. Note elements of the form $\overline{gH}$ with $%
g\in G$ are dense in $(\overline{G}/\overline{H})^{\symbol{94}}.$ Thus to
prove inversion in $(\overline{G}/\overline{H})^{\symbol{94}}$ is continuous
we need only check the case $f_{n}\rightarrow f,$ $f_{n}\in G,\overline{%
f_{n}H}\rightarrow \overline{fH},$ and $\overline{fpH}=\overline{fpH}=%
\overline{H}.$

We must show $\overline{f_{n}^{-1}H}\rightarrow \overline{pH}.$ Note $%
d(f,f_{n})=d(ff_{n}^{-1},id).$ Thus $\overline{ff_{n}^{-1}H}\rightarrow 
\overline{H}.$ Hence $\overline{f_{n}^{-1}H}=\overline{pfH}\ast \overline{%
f_{n}^{-1}H}=\overline{pff_{n}^{-1}H}\rightarrow \overline{pH}.$ Let $D$
denote the Hausdorf metric on $(\overline{G}/\overline{H})^{\symbol{94}}.$
Define as in Lemma \ref{compat} a metric $D^{\symbol{94}}$ on $(\overline{G}/%
\overline{H})^{\symbol{94}}$ via $D^{\symbol{94}}(\overline{fH},\overline{gH}%
)=D(\overline{fH},\overline{gH})+D((\overline{fH})^{-1},(\overline{gH}%
)^{-1}).$

Suppose $\overline{f_{n}H}$ is $D^{\symbol{94}}$ Cauchy. Let $\overline{%
f_{n}H}\rightarrow \overline{fH}$ and $\overline{(f_{n}H})^{-1}\rightarrow 
\overline{gH}$ with $\overline{\{fH},\overline{gH}\}\subset C_{H}(X,X).$
Then $\overline{f_{n}H}\ast (\overline{f_{n}H})^{-1}\rightarrow \overline{fH}%
\ast \overline{gH}=\overline{H}.$ Similarly $\overline{gH}\ast \overline{fH}=%
\overline{H}.$ Hence $\overline{\{fH},\overline{gH}\}\subset (\overline{G}/%
\overline{H})^{\symbol{94}}.$
\end{proof}

\begin{theorem}
\label{comgrp}The function $\phi :G/H\rightarrow $ $(\overline{G}/\overline{H%
})^{\symbol{94}}$ defined via $\phi (gH)=\overline{gH}$ is a homomorphism.
If $G$ is closed in $H(X,X)$ then $\phi $ is an epimorphism. If $H$ is
closed in $G$ then $\phi $ is a monomorphism and an embedding. Consequently
if $H$ is closed in $G$ and $G$ is closed in $H(X,X)$ then $\phi $ is a
homeomorphism and in particular $G/H$ is completely metrizable.
\end{theorem}

\begin{proof}
$\phi $ is well defined since by Corollary \ref{cor1} sets of the form $%
\overline{fH}$ determine a partition of $C(X,X).$ To show $\phi $ is a
homomorphism first note $\phi (g_{1}H)(g_{2}H)=\phi (g_{1}g_{2}H)=\overline{%
g_{1}g_{2}H}.$ However by Lemma \ref{th1} $\overline{g_{1}g_{2}H}=\ast (%
\overline{g_{1}H},\overline{g_{2}H}).$

Suppose $G$ is closed in $H(X,X)$ and $\overline{fH}\in (\overline{G}/%
\overline{H})^{\symbol{94}}.$ Let $\overline{pH}\in (\overline{G}/\overline{H%
})^{\symbol{94}}$ satisfy $\overline{pfH}=\overline{fpH}=\overline{H}.$
Choose a sequence $p_{n}\in \overline{pH}$ such that $fp_{n}\rightarrow id.$
Choose a sequence $f_{n}\in \overline{fH}.$ Such that $pf_{n}\rightarrow id.$
By Lemma \ref{conv} $f_{n}\rightarrow f.$ Fixing $\varepsilon >0$ since $%
p_{n}f_{n}\rightarrow id$ it follows that $f_{n}$ is eventually an $%
\varepsilon $ map.

Thus, those maps of $\overline{fH}$ which are $\varepsilon $ maps form an
open and dense subspace of $\overline{fH}.$ Since this holds for each $%
\varepsilon ,$ by the Baire category theorem there is a dense subset of
homeomorphisms within $\overline{fH}.$ In particular since $G=\overline{G}$
then $G\cap \overline{fH}\neq \emptyset .$ Choosing $g\in G\cap \overline{fH}
$ we have $\phi (gH)=\overline{gH}=\overline{fH}.$

Suppose $H$ is closed in $G.$ Then the closure of $gH$ within $G$ is $gH$
and therefore we apply Lemma \ref{lem2} to obtain a Hausdorf metric $d_{1}$
on $G/H$ inherited from the uniform metric of $G.$ By Cor\ref{cor1} we may
use a Hausdorf metric $d_{2}$ on $(\overline{G}/\overline{H})^{\symbol{94}}$
inherited from the uniform metric of $\overline{G}.$ Note $%
d_{1}(g_{1}H,g_{2}H)=d_{2}(\overline{g_{1}H},\overline{g_{2}H})=d_{2}(\phi
(g_{1}H),\phi (g_{2}H))$ Thus if $H$ is closed in $G$ then $\phi $ is an
embedding.
\end{proof}

\begin{corollary}
\label{bigcor}$M(H_{\infty })$ admits a compatible complete metric 
\begin{equation*}
D:M(H_{\infty })\times M(H_{\infty })\rightarrow M(H_{\infty })
\end{equation*}
determined by the condition $D(fH_{\infty }^{0},gH_{\infty
}^{0})<\varepsilon $ iff there exists $f_{1}\in fH_{\infty }^{0},f_{2}\in
f^{-1}H_{\infty }^{0},g_{1}\subset gH_{\infty }^{0}$ and $g_{2}\in
g^{-1}H_{\infty }^{0}$ such that $d(f_{1},g_{1})<\varepsilon $ $\ $and $%
d(f_{2},g_{2})<\varepsilon $. (Note $d(p,q)<\varepsilon $ iff $\left|
p(x)-q(x)\right| <\varepsilon $ $\forall x\in D^{2}).$
\end{corollary}

\begin{proof}
Observe by definition $H_{\infty }$ is closed in $H(D^{2},D^{2}).$ By Lemma 
\ref{close} $H_{\infty }^{0}=\cap H_{n}^{0}$ and hence $H_{\infty }^{0}$ is
closed in $H_{\infty }$ since each of $H_{n}^{0}$ is closed in $%
H(D^{2},D^{2}).$ Now apply Theorem \ref{comgrp} and harness the metric from
the proof of Corollary \ref{cor2}.
\end{proof}

\section{A left ordering on $M(H_{\infty })$ compatible with its topology}

An ordered group is \textbf{left ordered} if $f<g\Rightarrow hf<hg.$ In \cite
{f1},\cite{w1} Dehornoy's left ordering \cite{deh2} on $B_{n}$ is
interpreted geometrically. Letting $\gamma _{i}\subset D^{2}$ denote the
line segment connecting $z_{i}$ and $z_{i+1},$ the question of whether a
given braid $gH_{n}^{0}\neq H_{n}^{0}$ is positive or negative is settled at
the smallest index $i$ for which $\gamma _{i}\notin (gH_{n}^{0})(\gamma
_{i}).$ Consequently, (recalling the monomorphism $k_{n}:B_{n}%
\hookrightarrow B_{n+1}$) for each $n$ there exists a left ordering $<_{n}$%
on $B_{n}$ satisfying the following conditions:

\begin{enumerate}
\item  If $fH_{n}^{0}<_{n}gH_{n}^{0}$ then $%
k_{n}(fH_{n}^{0})<_{n+1}k_{n}(gH_{n}^{0})$

\item  If $fH_{n}^{0}<_{n}gH_{n}^{0}$ and $pH_{n}^{0}\in B_{n}$ then $%
pfH_{n}^{0}<_{n}pgH_{n}^{0}$

\item  If $fH_{n}^{0}\neq gH_{n}^{0}$ then the question of whether $%
fH_{n}^{0}<_{n}gH_{n}^{0}$ or $gH_{n}^{0}<_{n}fH_{n}^{0}$ is settled by
inspection of the arc isotopy classes $(fH_{n}^{0})(\alpha _{i})$ and $%
(gH_{n}^{0})(\alpha _{i})$ for the smallest index $i$ satisfying $%
(fH_{n}^{0})(\alpha _{i})\neq (gH_{n}^{0})(\alpha _{i}).$
\end{enumerate}

\begin{theorem}
\label{main}There exists a left ordering on $M(H_{\infty })$ such that the
order topology is compatible with the quotient topology.
\end{theorem}

\begin{proof}
Given $\{gH_{\infty }^{0},fH_{\infty }^{0}\}\subset M(H_{\infty })$ declare $%
gH_{\infty }^{0}<fH_{\infty }^{0}$ iff there exists $n$ such that $%
gH_{n}^{0}<_{n}fH_{n}^{0}.$ Note by Lemma \ref{close} given $\{fH_{\infty
}^{0},gH_{\infty }^{0}\}\subset M(H_{\infty })$ then $fH_{\infty
}^{0}=gH_{\infty }^{0}$ iff $fH_{n}^{0}=gH_{n}^{0}\forall n.$ It follows
easily that $<$ is a left ordering on $M(H_{\infty }^{0}).$ For
compatibility suppose $g\in (g_{1}H_{\infty }^{0},g_{2}H_{\infty }^{0}),$ an
open interval in the order topology. Choose $n$ such that $%
g_{1}H_{n}^{0}<gH_{n}^{0}<g_{2}H_{n}^{0}$ and $%
g_{2}^{-1}H_{n}^{0}<g^{-1}H_{n}^{0}<g_{1}^{-1}H_{n}^{0}.$ Choose $\delta >0$
such that for $f\in H_{n}$ if $d(f,g)<\delta $ and $d(f^{-1},g^{-1})<\delta $
then $fH_{n}^{0}=gH_{n}^{0}.$ Taking the metric $D$ on $M(H_{\infty }^{0})$
from Corollary \ref{bigcor}$,$ it follows that the open metric ball $%
B(gH_{\infty }^{0},\delta )\subset (g_{1}H_{\infty }^{0},g_{2}H_{\infty
}^{0}).$

Conversely suppose $U$ is open in $M(H_{\infty }).$ Let $qH_{\infty }^{0}\in
U.$

Let $V=\{q^{-1}hH_{\infty }^{0}|hH_{\infty }^{0}\in U\}.$ Note $V$ is open (
since $M(H_{\infty }^{0})$ is a topological group) and $H_{\infty }^{0}\in
V. $ Choose $\varepsilon >0$ such that $B_{D}(H_{\infty }^{0},\varepsilon
)\subset V.$

Choose $k$ such that diam$(D^{2}\backslash D_{k})<\varepsilon .$ Choose $%
fH_{\infty }^{0}\in B_{D}(H_{\infty }^{0},\varepsilon )$ such that $%
f_{D_{k}}=id_{D_{k}}.$

Suppose $f^{-1}H_{\infty }^{0}<gH_{\infty }^{0}<fH_{\infty }^{0}.$ Then for
each $n$ $f^{-1}H_{n}^{0}\leq gH_{n}^{0}\leq fH_{n}^{0}.$

Suppose $1\leq i\leq k$ and $n$ is any positive integer. Then 
\begin{equation*}
(fH_{n}^{0})(\alpha _{i})=(f^{-1}H_{n}^{0})(\alpha _{i})=H_{n}^{0}(\alpha
_{i}).
\end{equation*}
Hence $H_{n}^{0}(\alpha _{i})=gH_{n}^{0}(\alpha _{i})$ since otherwise we
contradict one of $f^{-1}H_{n}^{0}\leq gH_{n}^{0}$ or $gH_{n}^{0}\leq
fH_{n}^{0}.$

Thus $\forall n$ and each $i\leq k$ $\alpha _{i}\in (gH_{n}^{0})(\alpha
_{i}).$ By Corollary \ref{nicecor} select $g^{\symbol{94}}\in gH_{\infty
}^{0}$ such that $g_{D_{k}}^{\symbol{94}}=id_{D_{k}}.$ Hence $D(gH_{\infty
}^{0},H_{\infty }^{0})<\varepsilon .$

This shows the open interval $(f^{-1}H_{\infty }^{0},fH_{\infty
}^{0})\subset V.$ Hence (since left composition preserves the ordering) the
open interval $(qf^{-1}H_{\infty }^{0},qfH_{\infty }^{0})\subset U.$

Thus the order topology of $M(H_{\infty })$ is compatible with the quotient
topology.
\end{proof}

\end{document}